\documentclass{article}

\usepackage{graphicx}
\usepackage{epstopdf}
\usepackage{amsmath}
\usepackage{amssymb}
\usepackage{amsthm}
\usepackage{url}
\usepackage{lineno}

\theoremstyle{plain}
\newtheorem{theorem}{Theorem}
\newtheorem{lemma}[theorem]{Lemma}
\newtheorem{proposition}[theorem]{Proposition}

\newtheorem{corollary}[theorem]{Corollary}
\newtheorem{remark}[theorem]{Remark}
\newtheorem{question}[theorem]{Question}
\newtheorem*{definition}{Definition}
\newtheorem{conjecture}{Conjecture}

\newcommand{\Hth}{\mathbb{H}^3}

\newcommand{\SLTC}{SL(2, \mathbb{C})}
\newcommand{\BB}{\mathcal{B}}

\newcommand{\FF}{\mathcal{F}}
\newcommand{\HH}{\mathcal{H}}
\newcommand{\TT}{\mathcal{T}}

\newcommand{\nat}{\mathbb{N}}
\newcommand{\RR}{\mathbb{R}}
\newcommand{\Z}{\mathbb{Z}}

\newcommand{\TorusxI}{\mathbb{T}^2 \times I}

\newcommand{\SxD}{\mathbb{S}^1 \times D^2}
\newcommand{\DxS}{\SxD}

\newcommand{\KxxI}{K^2\tilde{\times}I}

\newcommand{\NP}{NP}
\newcommand{\coNP}{\MakeLowercase{co}NP}
\newcommand{\Mob}{M\"{o}bius strip}
\newcommand{\homeo}{\cong}

\DeclareMathOperator{\interior}{int}
\DeclareMathOperator{\image}{im}

\begin{document}

\title{On the complexity of cusped non-hyperbolicity}
\author{Robert~Haraway~III \and Neil~R.~Hoffman}
\date{\today}
\maketitle 

\begin{abstract}
  We show that the problem of showing that a cusped 3-manifold $M$ is not hyperbolic is in \NP{}, assuming {\sc $S^3$ recognition} is in \coNP{}.
  To this end, we show that {\sc Irreducible toroidal recognition} lies in \NP{}.
  Along the way we unconditionally recover {\sc Satellite knot recognition} lying in \NP{}.
  This was previously known only assuming the Generalized Riemann Hypothesis.
  Our key contribution is to certify closed essential normal surfaces as essential in polynomial time in compact orientable irreducible $\partial$-irreducible triangulations.
  Our work is made possible by recent work of Lackenby showing several basic decision problems in 3-manifold topology are in \NP{} or \coNP{}.
\end{abstract}

\section*{Declarations and Acknowledgements}
Both authors declare no conflict of interest.

Both authors thank Bus Jaco and Ben Burton for separate summaries of crushing;
William Petterssen for conversations at the start of this project;
and the MATRIX Institute who hosted these conversations.
They also thank the referee for careful reading and extensive constructive criticism.
Robert Haraway thanks Nathan Dunfield and Marc Lackenby for questions and helpful comments about surface bundle recognition.
Neil Hoffman thanks Gary Miller for questions about $S^3$ recognition that led to this project.
Neil Hoffman was also partially supported by Simons Foundation grant \#524123.

\section{Introduction}\label{intro}

\subsection{Our main results}
One useful property to know about a 3-manifold is whether or not it is hyperbolic---that is, whether or not it is homeomorphic to a finite volume quotient of $\mathbb{H}^3$ by a group of isometries acting properly discontinuously.
The decision problem of {\sc hyperbolicity} is computable (\cite{Kuperberg14,Manning02}), but not much is known about its complexity class.
We will focus on the class of cusped orientable 3-manifolds.
Here a \emph{cusped} 3-manifold is a compact orientable 3-manifold with nonempty boundary consisting of a finite disjoint union of tori. We stress that all 3-manifolds we discuss will be orientable here and throughout the paper.
More specifically, we seek to understand the complexity of certifying a negative answer to this decision problem.
The main result of this paper is as follows:

\begin{theorem}\label{thm:main}
  Suppose {\sc $S^3$ recognition} lies in \coNP{}.
  Then among cusped 3-manifolds, {\sc hyperbolicity} lies in \coNP{}.
\end{theorem}

We have taken care to minimize our result's dependence on a \coNP{}\ solution to {\sc $S^3$ recognition}, the decision problem to recognize $S^3$ among triangulations of manifolds.
Of course, this problem is known to be decidable (\cite{rubinstein1995algorithm,thompson1994thin}) and lie in \NP{}\ (\cite{Ivanov08,Schleimer11}), but the hypothesis is about certificates for affirming a negative result.
If we restrict to the set of irreducible manifolds, we may obtain the unconditional result:

\begin{theorem}
  Among irreducible orientable cusped 3-manifolds, {\sc hyperbolicity} lies in \coNP{}.
\end{theorem}

Along the way, we recover unconditionally the following theorem of Baldwin and Sivek \cite[Thm.~1.3]{BaldwinSivek17}, removing its dependence upon the Generalized Riemann Hypothesis.
We follow the convention that {\sc satellite knot recognition} is the decision problem which takes as its input a knot (here a set of $n$ edges in embedded in $S^3$) and determines if it has a toroidal complement.
\begin{theorem}\label{thm:satellite}
  {\sc Satellite knot recognition} lies in \NP{}.
\end{theorem}

This follows as a special case of the following theorem.
We say the {\sc irreducible toroidal recognition} problem is the decision problem to decide if an irreducible manifold is toroidal.

\begin{theorem}\label{main_thm_toroidal}
  {\sc Irreducible toroidal recognition} lies in \NP{}.

  In other words, if $\TT$ triangulates an orientable closed or cusped 3-manifold that is also irreducible and toroidal, then there is a certificate that $\TT$ is irreducible and toroidal that is verifiable in time polynomial in $|\TT|$.
\end{theorem}

We structure our nonhyperbolicity certificate according to a celebrated theorem of Thurston.
Assuming $M$ is both irreducible and neither a solid torus nor $\TorusxI$, Thurston \cite{Thurston82} showed that $M$ is either Seifert fibered, toroidal, or hyperbolic.
Typically, the theorem below is presented 
adhering to 
that trichotomy, but the following reformulation is more useful for disproving hyperbolicity.

\begin{theorem}[{\cite[Thm.~2.3]{Thurston82}}]\label{thm:thurston}
  Let $M$ be a cusped 3-manifold.
  $M$ is non-hyperbolic if and only if it admits an essential connected compact surface of nonnegative Euler characteristic.
\end{theorem}

There are but seven connected compact surfaces of nonnegative Euler characteristic: the projective plane $P^2$, the sphere $S^2$, the disk $D^2$, the Klein bottle $K^2$, the torus $T^2$, the \Mob\ $M^2$, and finally the annulus $A^2$.
Thus, to prove our main theorem, it will suffice to give, for any triangulation $\TT$ of a non-hyperbolic cusped 3-manifold, a certificate verifiable in time polynomial in $\TT$ (that is, polynomial in its number of tetrahedra) that $\TT$ admits such a surface.
The certificates we give are either normal such surfaces or appeals to the works of others, which themselves appeal in no small part to normal surfaces.
Thus we will begin the body of this paper with a review of normal surfaces, and of useful intermediate results from \cite{AHT06} and \cite{Lackenby16}.
Then, after building up the necessary intermediate certificates, we give the certificates advertised above.
Before doing so, however, we now give a brief review of the literature.

\subsection{Brief review of some relevant literature}
Haken's foundational result that {\sc unknot recognition} is decidable \cite{Haken} established the importance of algorithms in low dimensional topology using normal surface theory.
Hass, Lagarias and Pippenger \cite{HLP99} analyzed the complexity of this problem and showed that {\sc unknot recognition} lies in \NP{}.
Their paper contains a number of other foundational results---namely, the bounds on the size of fundamental and vertex normal surfaces' coordinates---which continue to be relied upon or imitated (for example, later in this paper).

We point the reader to \cite{JT95} for many of the original solutions to important decision problems utilizing normal surface theory, and we also point the reader to Matveev's book \cite{Matveev07} for further background.

In Theorem 3 of \cite{Ivanov08}, Ivanov showed, among many other things, that {\sc $\SxD$ recognition} lies in \NP{}.
We appeal to this result frequently.
This is distinct from {\sc unknot recognition}, since in the latter problem one is allowed the assumption during verification that the given 3-manifold is known to be irreducible.
One could show that {\sc $\SxD$ recognition} lies in \NP{}\ more or less directly using Schleimer's result in \cite{Schleimer11} that {\sc $S^3$ recognition} lies in \NP{}, after appealing to \cite{AHT06} as Lackenby did to take exteriors in polynomial time.

Reducibility is where we assume that {\sc $S^3$ Recognition} lies in \coNP{}.
This assumption follows from the Generalized Riemann Hypothesis by the work of Zentner in \cite{Zentner16}, using his own work on splicing irreducible representations into $\SLTC$, and using Kuperberg's work in \cite{Kuperberg14} on turning such representations into nonabelian representations into $SL(2,\mathbb{F}_p)$ for sufficiently small $p$; one bounds the size required on such a prime $p$ using the Generalized Riemann Hypothesis.

Lackenby in \cite{Lackenby16} has shown that among irreducible cusped 3-manifolds, {\sc $\SxD$ recognition} lies in \coNP{}.
In fact, much more is true---he showed that {\sc irreducible knotted recognition} lies in \NP{}.
That is, if $M$ is irreducible and not $\SxD$, then there is a certificate verifying \emph{both} of these properties simultaneously.

For their satellite knot certificate, in \cite{BaldwinSivek17} Baldwin and Sivek made an appeal to the work of Berge \cite{Berge91} and Gabai \cite{Gabai89,Gabai90} on classifying knots in solid tori.
However, Baldwin and Sivek got a certificate using representations into $\SLTC$.
Thus, their certificate's size bound depends on Kuperberg's work, and hence on the Generalized Riemann Hypothesis.
Our appeals to normal surface theory removes this dependence in a number of relevant cases.
Nevertheless, we still depend on Berge and Gabai.

\section{Background}

To provide the minimal background needed in complexity theory, the reader should consult Schleimer's brief summary in \cite{Schleimer11}.
For our purposes, a decision problem \emph{lies in \NP{}} if there is a polynomial time verifiable certificate, or proof, of an affirmative answer.
Likewise, a decision problem \emph{lies in \coNP{}} if there is a polynomial time verifiable certificate of a negative answer.

When constructing such a certificate, one can assume more about a given input than one can at verification time.
A trivial general example is that during construction, one can assume the answer to the decision problem is Yes (if the problem is in \NP).
But during verification, this is just what the certificate asserts.
If we were allowed to assume the answer were Yes during verification, we could dispense with the certificate entirely.
In general, during verification, the only assumptions you can make about the input are those derivable in polynomial time from the input itself and from the certificate.
We gave an example above: the distinction between {\sc unknot recognition} and {\sc $\SxD$ recognition}.
With {\sc unknot recognition}, knot exteriors in $S^3$ are all irreducible, so one may assume irreducibility during verification.
The same is not true for {\sc $\SxD$ recognition}; the input is any triangulation, and arbitrary triangulations might be reducible.
To assume irreduciblity during verification, one would have to include a short proof of it in the certificate.

We will work throughout in the PL category, and we will use the standard kind of triangulation, defined as follows:

\begin{definition}
  Given a disjoint union $\overline{\TT}$ of closed tetrahedra; and given a partition $\Pi$ of some subset of the faces of these tetrahedra into pairs; and, finally, for each pair $(f,f')$ in $\Pi$, given a cell-isomorphism or \emph{gluing} $\phi_{(f,f')}$ between $f$ and $f'$, one may construct a cell-complex by identifying the tetrahedra along the gluings, and remembering all this information.
  A cell-complex so constructed we call a \emph{triangulation}.
  All triangulations are constructed this way.
  The \emph{underlying space} of a triangulation is the topological space resulting from the identification.

  A triangulation $\TT$ with finitely many tetrahedra is a \emph{(material) triangulation} of a compact 3-manifold $M$ when $\TT$ is homeomorphic to $M$.
  This will be the convention throughout this paper and we will use $|\TT|$ to denote the number tetrahedra in the triangulation $\TT$.
  
  Finally, we use $\TT^{(0)}$,$\TT^{(1)}$,$\TT^{(2)}$ to denote the 0-,1-
  and 2-skeletons of $\TT$.
\end{definition}

Because $\TT$ is a representative of the manifold $M$ and is indeed the input to our algorithms, we will use $\TT$ to denote the manifold of interest throughout the paper.

\begin{remark}
  The above definition is more permissive than the definition of a simplicial complex, since it does not require cell-inclusion maps to be embeddings.
  For instance, there is a triangulation of the 3-sphere with one tetrahedron.
  But the only simplicial complex with one tetrahedron is the tetrahedron itself.
\end{remark}

\begin{remark}
  Another notion of triangulation is \emph{ideal triangulations}, which typically use fewer tetrahedra.
  As useful as ideal triangulations are in other contexts, we eschew them here.
  We mention ``material'' triangulations only to emphasize this.
  It is worth pointing out that, given an ideal triangulation $\TT$ of a 3-manifold $M$, one may associate a material triangulation $\TT'$ to $M$ in time polynomial in $\TT$ by taking a second barycentric subdivision and removing the tetrahedra around the ideal vertices of $\TT$.
  A construction due to Weeks (see \cite[pp.~469--470]{WeeksKT}) allows the converse operation, going from a material triangulation to an ideal triangulation.
\end{remark}

\begin{remark}
  Lackenby prefers to use handle structures instead of triangulations (for good reason; see \cite[\S 1.3]{Lackenby16}), but we relied mainly on triangulations in this work.
  These two kinds of structure are polynomially equivalent.
  That is, given a handle structure $H$ on a 3-manifold with a given total number of cells $h$, one may calculate in time polynomial in $h$ a triangulation $T$ of the same 3-manifold, and vice versa.
  Lackenby gives one such equivalence in \cite{Lackenby16}.
  Triangulations admitting splitting surfaces as in \cite[Chapter 4]{BurtonPhD} provide another such equivalence.
\end{remark}

Our work depends on the theory of normal surfaces.

\begin{definition}
  Suppose $\TT$ is an orientable triangulation.
  All the following definitions depend upon this choice of $\TT$.

   A \emph{normal isotopy} of a triangulation is an isotopy $\phi:\TT\times I \to \TT$ leaving the triangulation invariant---i.e.~such that, for every $t \in I$, the homeomorphism $x \mapsto \phi(x, t)$ is a cell-isomorphism from $\TT$ to itself.

  A \emph{normal disk} in a tetrahedron $\tau$ of $\TT$ is a properly embedded disk in $\tau$ normally isotopic to the intersection of $\tau$ with an affine plane transverse to $\tau$.
  If it separates one vertex from the others, or \emph{links} the vertex, then the disk is a triangle.
  Otherwise the plane separates two pairs of vertices, and the disk is a quadrilateral, or \emph{quad}.

  A \emph{normal surface in $\TT$} is a properly embedded surface $\Sigma \hookrightarrow \TT$ transverse to $\TT^{(2)}$ such that $\Sigma$ is the union of finitely many normal disks.

  A \emph{link} of a vertex $v$ of $\TT$ is a normal surface in $\TT$ consisting of one of each normal triangle disk linking $v$.
  These are ``trivial'' normal surfaces.

  The \emph{normal disk set of $\TT$} is the set $\Delta$ of normal isotopy classes of normal disks in $\TT$.
  It has cardinality $7\cdot |\TT|$.

  The \emph{(normal) coordinates} of a normal surface $\Sigma$ are constituted by the function $x_{\Sigma}: \Delta \to \nat$ defined by letting $x_{\Sigma}(d)$ be the number of normal disks of $\Sigma$ in the normal isotopy class $d$.

  The \emph{(total) weight} of a PL surface $\Sigma$ in $\TT$ transverse to $\TT^{(2)}$ is $|\Sigma \cap \TT^{(1)}|$.
\end{definition}

A normal surface is determined up to normal isotopy by its coordinates.
Therefore, one may combinatorially represent a normal surface via its coordinates.
That is to say, letting $\mathcal{N}(\TT)$ be the set of normal isotopy classes of normal surfaces in $\TT$, the function $x_\cdot: \mathcal{N}(\TT) \to \mathbb{N}^{\Delta(\TT)}$ taking such a class to the coordinates of one of its members is not only well-defined, but also is injective.
However, it is not surjective.
For instance, no normal surface $\Sigma$ has $x_\Sigma(q) > 0$ and $x_\Sigma(q') > 0$ for two distinct normal isotopy types $q,q'$ in a single tetrahedron of $\TT$.

Not every surface in $\TT$ is isotopic to a normal surface.
However, many surfaces important for 3-manifold topology can be so isotoped.
We will also regard those properties $\TT$ (or $M$) that can be computed from $\TT$ in polynomial time as inherent to $\TT$.
For example, if $\TT$ can be distinguished from $\TorusxI$ by a polynomial-time computable invariant such as number of cusps or homology, we say it can be distinguished via the empty certificate.
Our convention here would be to avoid unnecessary storage of certificates in such cases.

We refer the reader to \cite{JR03} for a simple ``shrinking'' normalization procedure that produces normal such surfaces from non-normal ones.

The set of normal surfaces is infinite, even when considered only up to normal isotopy.
However, normal surfaces can sometimes be added together.
Normal surfaces that cannot be expressed as such a sum end up forming a finite set.
That is the class of surfaces we will focus on.
While summation can be defined strictly with respect to vectors, we use geometric sum in Lemma \ref{lem:nonsepAx}.
Geometric sum uses the following definition.

\begin{definition}
  Two quads $q,q'$ in the same tetrahedron with different normal isotopy type are said to \emph{clash}.
  A pair of normal surfaces $X,Y$ is \emph{admissible} when no two normal quads from $X$ and $Y$ clash.
\end{definition}

If $S,S'$ are an admissible pair of normal surfaces, then $x_S + x_{S'}$ is the vector of a normal surface.
Briefly, this is because both $x_S,x_{S'}$ are solutions to a homogeneous linear system of equations, Haken's \emph{matching equations}; by homogeneity, $x_S + x_{S'}$ is therefore also a solution; and by admissibility of the pair, the formal solution is the vector of an actual normal surface.

Haken gave a geometric construction of this surface, which construction is called the geometric or Haken \emph{sum} $S + S'$.
Briefly, we may normally isotope $S$ so that $S,S'$ are transverse.
Suppose $\kappa$ is a component of $S \cap S'$.
Let $N$ be a small regular neighborhood of $\kappa$ in $\TT$.
We can establish a product structure on $N$ such that $N = \kappa \times D^2$ and such that for all $p \in \kappa$, $(S \cup S') \cap (\{p\} \times D^2)$ is $\times$ in the disc, i.e.~two properly embedded arcs intersecting once in their interiors.
We can replace $(S \cup S') \cap N$ with one of two surfaces, whose intersection with every such $\{p\} \times D^2$ is instead a pair of disjoint arcs with the same four boundary points as the $\times$.
This \emph{resolves} the intersection.
Haken showed that when $X,Y$ are admissible, exactly one of the two resolutions produces an immersed normal surface, i.e.~a (possibly non-disjoint) union of normal discs.
On performing all such resolutions, one gets an immersed normal surface with no self-intersections, i.e.~just a normal surface.
Any two such normal surfaces gotten by completely resolving the intersections are normally isotopic.
So $X+Y$, the geometric sum of $X$ and $Y$, is well-defined up to normal isotopy.

For Lemma \ref{lem:nonsepAx} we note here that if $[X]$ and $[Y]$ are the classes of $X$ and $Y$ in $H_2(\TT,\partial \TT; \mathbb{Z}/2 \mathbb{Z})$, then $[X+Y] = [X]+[Y]$, because resolutions don't change $\mathbb{Z}/2 \mathbb{Z}$ homology.
Resolutions also respect weight, so that we have $wt(X+Y) = wt(X) + wt(Y)$.
Finally, $X+Y$ has as many $d$-cells as $X$ and $Y$ combined, for all $d = 0,1,2$.
So in particular, $\chi(X+Y) = \chi(X) + \chi(Y)$.

Now we define the finite class of normal surfaces we use.

\begin{definition}
  A normal surface $\Sigma$ is \emph{fundamental} when, for any two normal surfaces $X,Y$,
  we have $X + Y = \Sigma$ if and only if $\{X,Y\} = \{\Sigma, \emptyset\}$.
\end{definition}

A common refrain in
 normal surface theory is that if there is an ``interesting'' surface, then there is a fundamental one.
We first give some examples of ``interesting'' surfaces.
We follow that by a proposition that fits the trope.
It is of fundamental importance to this paper, and is how we structure our certificate of non-hyperbolicity and its verification.
Its proof mostly just points the reader to \cite{Matveev07}.

\begin{definition}
  Suppose $\TT$ triangulates a compact 3-manifold.
  Suppose $\Sigma$ is a tame, connected, properly embedded surface in $\TT$.

  $\Sigma$ is an \emph{essential sphere} when it is an sphere that does not bound a ball.

  $\TT$ is \emph{irreducible} when it admits no essential sphere.

  $\Sigma$ is a \emph{compressing disk} when it is a disk whose boundary is essential in $\partial \TT$.

  $\TT$ is \emph{$\partial$-irreducible} when it admits no compressing disk.

  $\TT$ is \emph{0-efficient} when it admits no nontrivial normal surfaces of positive Euler characteristic.

  $\Sigma$ is \emph{compressible} when there is an essential curve $\gamma$ in $\Sigma$ and a tame embedded disk $D$ in $\TT$ such that $D\cap \Sigma = \gamma = \partial D$.

  $\Sigma$ is \emph{$\partial$-compressible} when there is an essential arc $\alpha$ in $\Sigma$ and a tame embedded disk $D$ in $\TT$ and an arc $\beta$ of $D$ such that $\alpha = D \cap \Sigma$, $\beta = D \cap \partial \TT$, $\alpha \cap \beta = \partial \alpha = \partial \beta$, and $\alpha \cup \beta = \partial D$.

  $\Sigma$ is \emph{$\partial$-parallel}
  when it is isotopic relative to its boundary to a subsurface of $\partial \TT$.

  $\Sigma$ is \emph{essential} (assuming it is not a sphere) when it is incompressible, $\partial$-incompressible, and not $\partial$-parallel.

  $\TT$ is \emph{toroidal} when it is irreducible and $\partial$-irreducible and admits an essential torus.
\end{definition}

\begin{proposition}\label{prp:fundFaults}

  Suppose $\TT$ is an orientable material triangulation of a 3-manifold.
  Let $\FF$ be the set of fundamental 
  normal surfaces in $\TT$ up to normal isotopy. Then
  \begin{itemize}
  \item if there is an embedded $P^2$, then there is one in $\FF$;
  \item otherwise, if there is an essential $S^2$, then there is one in $\FF$;
  \item otherwise, if there is a compressing $D^2$, then there is one in $\FF$;
  \item otherwise, if there is an essential $K^2$, then there is one in $\FF$;
  \item otherwise, if there is an essential $T^2$, then there is one in $\FF$;
  \item otherwise, if there is an essential $M^2$, then there is one in $\FF$;
  \item otherwise, if there is an essential $A^2$, then there is one in $\FF$.
  \end{itemize}
\end{proposition}

\begin{proof}
  Suppose $\TT$ is an orientable material triangulation.
  Let $\mathcal{F}$ be the set of fundamental normal surfaces (up to normal isotopy).

  \begin{description}

  \item[$P^2$:] If $\TT$ admits an embedded $P^2$, then by \cite[Theorem 4.1.12]{Matveev07}, $\TT$ admits a fundamental such surface.

  \item[$S^2$:] Otherwise, if $\TT$ admits an essential sphere, then again by \cite[Theorem 4.1.12]{Matveev07}, $\TT$ admits a fundamental such sphere, since $\TT$ admits no embedded $P^2$.

  \item[$D^2$:] Otherwise, if $\TT$ admits a compressing disk, then by the proof of \cite[Theorem 4.1.13]{Matveev07}, $\TT$ admits a fundamental such disk, since $\TT$ is irreducible.

  \item[$K^2$:] Otherwise, if $\TT$ admits an essential $K^2$, then by irreducibility there is a normal surface $\Sigma$ isotopic to this Klein bottle, e.g.\,\,by applying the shrinking procedure of \cite{JR03}.
    Suppose $\Sigma = S + S'$.
    We follow an argument similar to the proof of \cite[Lemma 6.4.7]{Matveev07}.
    Since $\TT$ is irreducible and orientable, and since $\Sigma$ is incompressible, by \cite[Lemma 3.3.30]{Matveev07} and \cite[Theorem 4.1.36]{Matveev07}, both $S$ and $S'$ are connected, incompressible, and not $S^2$ or $P^2$.
    Hence they both have Euler characteristic 0.
    Also, they are closed.
    Since $\Sigma$ is not orientable, at least one of $S,S'$ is not orientable.
    Without loss of generality, $S$ is not orientable; since it is closed and $\chi(S) = 0$, $S$ is a Klein bottle.
    If $S'$ is not empty, then the weight $w(S) = |S \cap \TT^{(1)}|$ is less than $w(\Sigma) = |\Sigma\cap\TT^{(1)}|$.
    Hence by descent on weight, one arrives at an embedded Klein bottle $\Sigma$ such that if $\Sigma = S + S'$, then one of $S$ or $S'$ is empty---that is, we arrive at a fundamental normal essential Klein bottle.

  \item[$T^2$:] Otherwise, suppose $\TT$ admits an essential $T^2$.
    By \cite[Lemma 6.4.7]{Matveev07}, $\TT$ admits either a fundamental normal such surface, or a fundamental normal Klein bottle $K$ with $K + K$ being an essential $T^2$.
    But this latter case would mean $K$ was an essential Klein bottle, contrary to assumption.
    Thus $\TT$ admits a fundamental normal essential $T^2$.

  \item[$M^2$:] Otherwise, if $\TT$ admits an essential \Mob\ $M^2$, then it admits a normal essential band $\Sigma$, again by \cite{JR03}.
    Just as before, if $\Sigma = F + F'$, then by the same lemmata, $F$ and $F'$ are connected, incompressible, $\partial$-incompressible, and not $S^2$, $P^2$, or $D^2$.
    Hence their Euler characteristics are 0.
    One of them, say $F$, has exactly one boundary component, and the other $F'$ is closed.
    If $F'$ were nonorientable, then it would be an embedded Klein bottle, contrary to assumption.
    Hence $F$ is nonorientable, and is thus an embedded M\"{o}bius band; and by the same lemmata as above is essential, being incompressible and $\partial$-incompressible.
    By descent as above, $\TT$ admits a fundamental normal embedded M\"{o}bius band.

  \item[$A^2$:] Otherwise, if $\TT$ admits an essential annulus, then it admits a normal such annulus $\Sigma$.
    By \cite[Lemma 6.4.8]{Matveev07}, there is a fundamental such annulus, since $\TT$ admits no M\"{o}bius band.
  \end{description}
  \qed
\end{proof}

The following bound of Hass, Lagarias, and Pippenger improves on the mere computability of $\FF$.

\begin{lemma}[Hass, Lagarias, Pippenger {\cite[{Lemma 6.1(2)}]{HLP99}}]\label{lem:hlp}
  If $\TT$ is a triangulation of a compact 3-manifold with $t$ tetrahedra, and $\FF$ is the set of coordinates of fundamental normal surfaces of $\TT$, then for every $f \in \FF$,
  \[ \max_{d \in \Delta(\TT)} f(d) \leq t\cdot 2^{7\cdot t + 2}.\]
\end{lemma}

The point about this bound is that the surfaces from the previous proposition admit representations of (bit-)size polynomial in $t$, since to represent a surface one may represent its coordinates.
One may represent these coordinates in place-value notation, which only requires space proportional to the logarithm of these coordinates.

The fact that one may represent such a normal surface with a polynomial amount of data suggests, together with Kneser-Haken finiteness, that the topological classification of a normal surface $F$ of total weight $W$ ought to be at least representable with an amount of data polynomial in $\TT$ and $\log W$.
In fact, much more is true.
We have the following remarkable corollary of Agol, Hass, and Thurston in \cite[Corollary 17]{AHT06}, another result to which we appeal frequently:

\begin{proposition}[Agol, Hass, Thurston]\label{prp:aht}
  Suppose $\TT$ is a triangulation with $t$ tetrahedra, and suppose $F$ is a normal surface in $\TT$ of total weight $W = |F \cap \TT^{(1)}|$.

  There is an algorithm that, in time polynomial in $\TT$ and $\log W$, determines the coordinate vectors of the normal isotopy types of components of $F$, determines the homeomorphism class of each such type, and determines how many of each such type there are.
\end{proposition}

\begin{figure}
  \includegraphics[width=\textwidth]{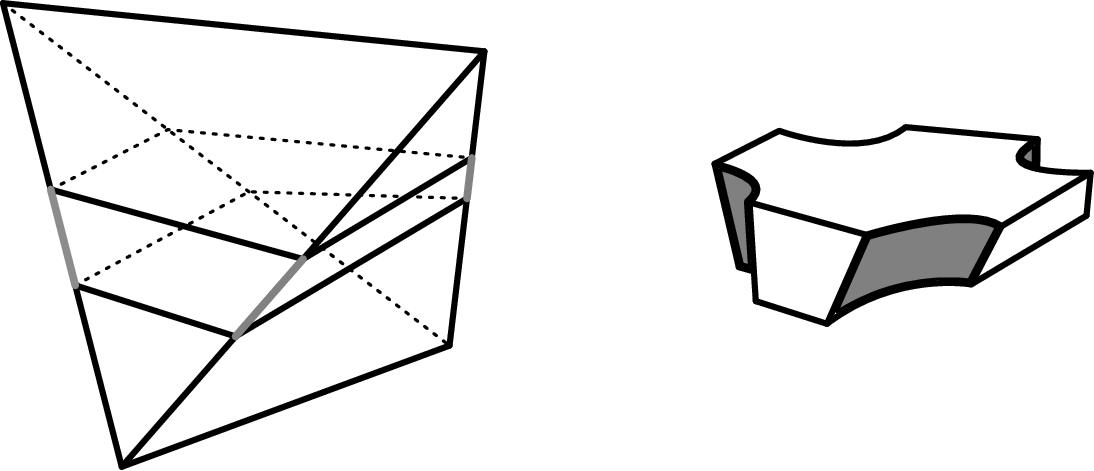}
  \caption{A parallelity 0-handle coming from a type IV cell}\label{fig:0_handle}
\end{figure}

\begin{figure}
  \includegraphics[width=\textwidth]{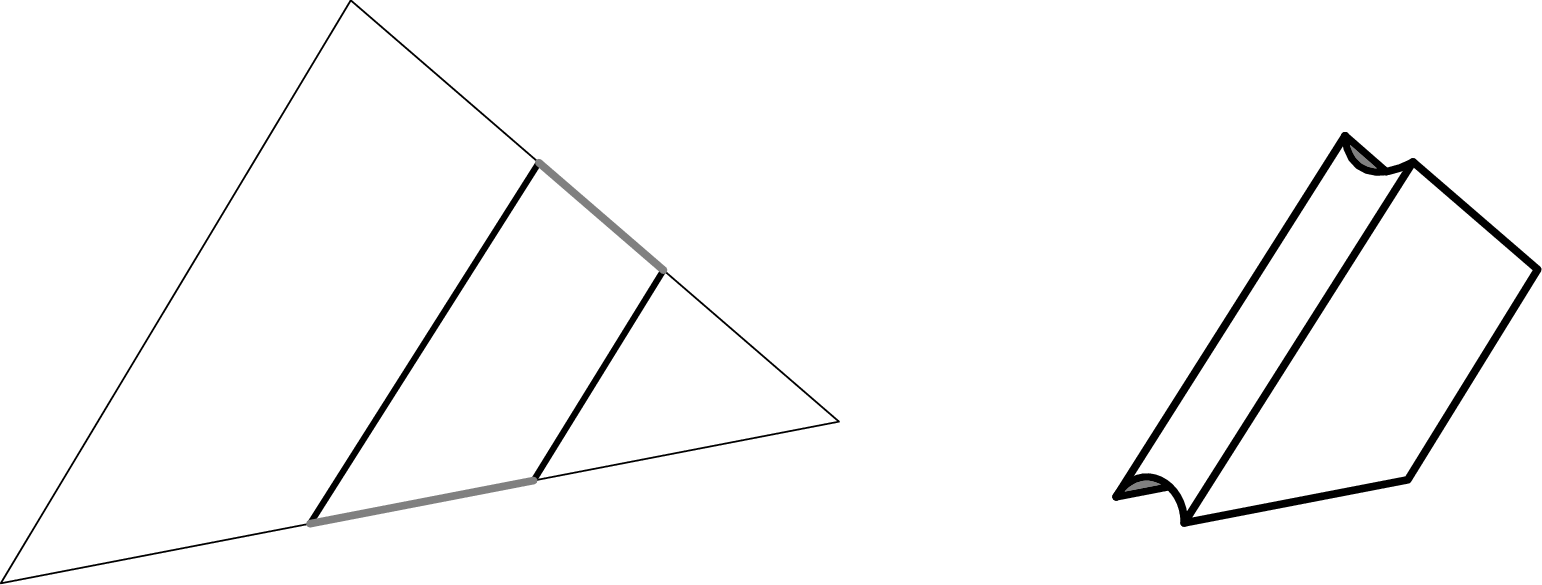}
  \caption{A parallelity 1-handle coming from a trapezoid}\label{fig:1_handle}
\end{figure}

\begin{figure}
  \includegraphics[width=\textwidth]{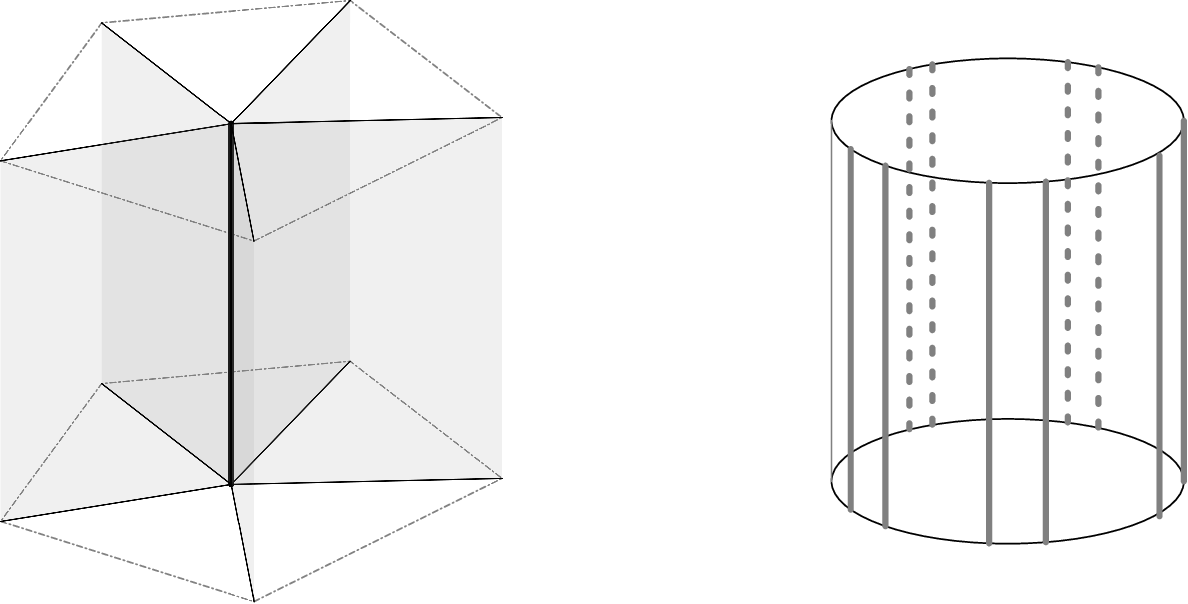}
  \caption{A regular neighborhood of an interstitial arc, and the parallelity 2-handle coming from the arc}\label{fig:2_handle}
\end{figure}

We also recall here a crucial intermediate result, a simple consequence of \cite[Theorem 9.3]{Lackenby16}, which Baldwin and Sivek also relied upon in \cite{BaldwinSivek17}.
The result relies on the notion of \emph{parallelity bundle}.
In general this is a notion to do with handle structures.
For us the particular kind of handle structure that motivates the definition is this.
Suppose $\Sigma$ is a normal surface in a triangulation $\TT$.
Let $C$ be the natural cellulation on the exterior of $\Sigma$ in $\TT$.
Let $\HH$ be the \emph{inverted} corresponding handle structure on the underlying space of $C$.
That is, let the 3-handles $\HH^3$ of $\HH$ be a small regular neighborhood of the vertices of $C$;
and for $0 \leq i < 3$ let the $(2-i)$-handles $\HH^{2-i}$ of $\HH$ be a small regular neighborhood of $C^{(i+1)} \setminus \interior \HH^{3-i}$ in $C$.
So 2-handles come from edges of $C$, 1-handles from faces, and 0-handles from 3-cells.
See Figures \ref{fig:0_handle}, \ref{fig:1_handle}, and \ref{fig:2_handle}.

Let $\FF^0$ be $\HH^1 \cap \partial \HH^0$ and $\FF^1$ be $\HH^2 \cap \partial \HH^0$.
These are ``islands'' and ``bridges'' in Matveev's terminology \cite[\S3.4]{Matveev07}.
The choice of the letter $\FF$ has nothing to do with fundamental surfaces, either.

With these conventions we make the following definition.
\begin{definition}[Parallelity bundle]
  Quoting from \cite[\S9.2]{Lackenby16}, a \emph{parallelity handle} $H$ of $\HH$ is a handle admitting ``a product structure $D^2 \times I$ such that
  \begin{enumerate}
  \item $D^2 \times \partial I = H \cap \Sigma$; [and]
  \item each component of $\FF^0 \cap H$ and $\FF^1 \cap H$ is $\beta \times I$ for a subset $\beta$ of $\partial D^2$.''
  \end{enumerate}
\end{definition}

\begin{theorem}[Lackenby]\label{thm:lknby_complement}
  There is an algorithm that takes, as its input, 
  \begin{enumerate}
  \item[(i)] a triangulation $\TT$, with $t$ tetrahedra, for a compact orientable manifold $M$; and
  \item[(ii)] a vector $\vec{S}$ for an orientable normal surface $S$ with no two components normally isotopic,
  \end{enumerate}
  and provides, as its output, the following data, in time that is bounded by a polynomial in $t \cdot \log |S \cap \TT^{(1)}|$.

  If $M'$ is the manifold that results from decomposing along $S$, and $S'$ is the two copies of $S$ in $\partial M'$, and $\BB$ is the parallelity bundle for the pair $(M',S')$ with its induced handle structure, then the algorithm determines the following information: a handle structure for $cl(M' - \BB)$ and, for each component $B$ of $\BB$, 
  \begin{enumerate}
  \item[(i)] the genus and number of boundary components of its base surface;
  \item[(ii)] whether $B$ is a product or twisted $I$-bundle; and
  \item[(iii)] the location of $\partial_v B$ in $cl(M'-\BB)$.
  \end{enumerate}
\end{theorem}

With the above theorem we can state the following proposition. We placed the proof of the proposition in an appendix. 

\begin{proposition}\label{prp:exterior}
  There is an algorithm that takes as its input both a compact orientable connected triangulation $\TT$ of $t$ tetrahedra and a connected normal surface $S$ in $\TT$ given as a vector, and provides as its output a triangulation of an exterior of $S$, and that moreover provides this output in time polynomial in $t$ and $|\chi(S)|$.
\end{proposition}
The example of a large number of parallel normal discs in a tetrahedron shows it is necessary for $S$ to be connected.

We provide a proof of Proposition \ref{prp:exterior} in an appendix, since it follows fairly easily from Lackenby's work.

The above operation due to Lackenby can triangulate the exterior of, say, a normal sphere of weight $2^n$ with $O(|\TT|+n)$ tetrahedra, where $\TT$ is the original triangulation.
As good as this is, it may not reduce the number of tetrahedra.
On the other hand, the \emph{crushing} operation always produces a smaller triangulation. 
This procedure was introduced by Jaco and Rubinstein in
\cite{JR03} and refined by Burton in \cite{BurtonCrush}. The implementation of crushing below is in some sense a blend of the two. 
We omit some of Jaco and Rubinstein's restrictive conditions to set up a crushing procedure. This reduces the process to relatively easy-to-implement operations at the expense of possibly over-simplifying the prime
decomposition of the manifold. In that sense, the result of crushing may not account for all components in the prime decomposition of the exterior.
However, when crushing along spheres or discs, the differences from the exterior are well-understood and will be discussed below.
Furthermore, we will only use this operation to certify homeomorphism between a triangulation and a 0-efficient triangulation.

\newcommand{\cut}{\downharpoonright}
Before defining crushing, we start with something more straightforward.
Suppose $\Sigma$ is a normal surface in the closed triangulation $\TT$.
Let $\Sigma^\ast$ be $\Sigma$ with a disjoint union of vertex links, so that $\Sigma^\ast$ contains one of each such vertex link as a component.
Then the exterior $\cut\Sigma^\ast$ of $\Sigma^\ast$ in $\TT$ naturally inherits a cell structure from $\TT$ with \emph{interstitial} cells of five types.
The type I cells are truncated tetrahedra.
The type II cells are truncated half-tetrahedra, or purses.
The type III and IV cells are triangular and quadrilateral prisms, respectively.
Finally, the type V cells are tetrahedra incident to vertices of $\TT$.
The interesting components of $\cut \Sigma^\ast$ are those without type V cells.
So we want to ignore cells of type V.
\newcommand{\crack}{\$}
\begin{definition}
  If $\Sigma$ is normal in $\TT$, and $\Sigma^\ast$ is $\Sigma$ with enough vertex links added to have one of every vertex link, then $\cut\Sigma^\ast$ with all type V cells omitted is $\TT$ \emph{cracked along} $\Sigma$, for which we write $\crack\Sigma$.
\end{definition}

Now, to crush, from $\crack \Sigma$ we first throw out the cells of types III and IV in $\$\Sigma$.
Next, suppose $c$ is a type II cell of $\$\Sigma$.
Then $c$ has two hexagonal faces $h_0, h_1$ incident along an edge $e$.
Since $\TT$ is closed, some hexagonal face $h_0'$ of a cell of $\$\Sigma$ different from $h_0$ is glued to $h_0$.
If $h_0' = h_1$, then we throw out $c$ from $\$\Sigma$.
Otherwise, let $\phi_0$ be the cell-isomorphism pairing $h_0'$ to $h_0$ for $\$\Sigma$.
Likewise, let $\phi_1$ be the cell-isomorphism similarly pairing $h_1$ to $h_1'$ for $\$\Sigma$.
Finally, let $\phi_c$ be the cell-isomorphism from $h_0$ to $h_1$ that ``folds along'' their common edge $e$.
We throw out $c$, and replace $\phi_0, \phi_1$ with the single face-pairing $\phi_1 \circ \phi_c \circ \phi_0$.
Throwing out type II cells this way, we eventually get a cellulation with only type I cells.
At this point, we just cone off the triangulated boundary to get a triangulation.
\newcommand{\crush}{!}
\begin{definition}
  This final coned-off triangulation is \emph{$\TT$ crushed along $\Sigma$}, which we write as $\crush \Sigma$.
\end{definition}

Crushing along a sphere or disc at most can do the following to a triangulation see \cite[Cor.~5]{BurtonCrush}:
\begin{theorem}[Burton]\label{thm:burton_crush}
  Suppose $\TT$ is a 3-triangulation containing no two-sided $P^2$.
  Let $S$ be a normal sphere or disc in $\TT$.
  Let $\TT_{JR}$ be $\TT$ crushed along $S$.
  Then $\TT_{JR}$ is a valid 3-triangulation obtained from $\TT$ by zero or more of the following operations:
  \begin{itemize}
  \item undoing connected sums, i.e.~surgering along spheres;
  \item cutting open along properly embedded discs;
  \item filling boundary spheres with 3-balls;
  \item deleting components homeomorphic to $D^3$, $S^3$, $RP^3$, $L(3,1)$, $S^2\times S^1$, or $S^2\tilde{\times}S^1$.
  \end{itemize}
\end{theorem}

\section{Certificates}

\subsection{Three-sphere certificates}
Rubinstein and Thompson both describe algorithms for 3-sphere recognition \cite{rubinstein1995algorithm,thompson1994thin}.
Later, it was shown that such algorithms could also produce an appropriately sized certificate affirming a manifold is indeed $S^3$.

\begin{theorem}[Ivanov \cite{Ivanov08}, Schleimer \cite{Schleimer11}]\label{thm:s3np}
  {\sc $S^3$ recognition} lies in \NP{}.
\end{theorem}

The companion question remains open.

\begin{conjecture}\label{conj:s3conp}
  {\sc $S^3$ recognition} lies in \coNP{}.
\end{conjecture}

Zentner \cite[Theorem 11.2]{Zentner16}
proved the following, providing strong evidence for the conjecture.

\begin{theorem}[Zentner]
  {\sc $S^3$ recognition} lies in \coNP{}, provided that the Generalized Riemann Hypothesis holds.
\end{theorem}

\subsection{Isomorphism signatures}\label{subsec:isosigs}

It may happen on occasion that we will want our triangulations to have a particular form---for instance, we may want to have a triangulation that induces a one-vertex triangulation on a boundary component.
One may simplify a triangulation to have such properties in polynomial time, and then one may generate certificates for the simplified triangulation.
However, these are not certificates for the original triangulation.

To promote these certificates to certificates of the original triangulation, it suffices to give a polynomial-sized certificate that the simplified triangulation triangulates the same underlying manifold as does the old triangulation.
It is almost obvious that one may do this for polynomial-length sequences of Pachner moves, layerings, close-the-book moves, and other ``atomic'' modifications of triangulations.
However, it is not clear how to name such a sequence in a way that is invariant under isomorphism of the triangulation.
Moreover, a representation of, say, a normal surface in a triangulation as a vector of numbers requires some choice of ordering on the normal disk-types in that triangulation.
One can consistently specify a natural such ordering given an ordering of the tetrahedra of the triangulation, and given, for each tetrahedron, an ordering of its vertices.
But then to use that vector as a certificate, one must make sure to make that particular collection of choices of ordering for one's triangulation.

Therefore, following \cite[\S 3]{Schleimer11}, let us say that a triangulation together with such choices is a \emph{labelled triangulation}.
Let two such structures be \emph{equivalent} when there is an isomorphism between the triangulations preserving the orderings (or \emph{labellings}).
In \cite{Burton11} Burton has constructed a injective \emph{signature} function from the set of equivalence classes of labelled triangulations to the set of bit-sequences; and has constructed a function, the \emph{isomorphism signature} function, from the class of labelled triangulations to the set of bit-sequences, whose level classes are in bijection with \emph{combinatorial isomorphism classes of triangulations}; and, most importantly for us, these functions may be computed in polynomial time in their inputs.
This gives canonical ways of putting coordinates on normal surfaces and other such objects in triangulations.

Therefore, throughout the work below, we assume that every triangulation is labelled, and in fact has the labelling whose signature is its isomorphism signature, which is just the lexicographic minimum of the set of signatures over all its possible labellings---we assume it has its \emph{canonical labelling}.
When implementing the algorithms below, the first thing one should do with a newly constructed triangulation or with an input triangulation is to endow it with its canonical labelling.

\subsection{Simplifying triangulations}\label{subsec:simpl}

A standard ``close-the-book-and-layer'' algorithm (e.g.\,\,the algorithm in \cite{SnapPy}, \texttt{close\textunderscore cusps.c}) returns, in time polynomial in $\TT$, a new triangulation $\TT'$ such that $\TT'$ has at most three more tetrahedra in each boundary component
than $\TT$ and such that $\TT'$ induces one-vertex triangulations on $\partial \TT'$; and also returns a proof $P$ of size polynomial in $\TT$ that $\TT$ and $\TT'$ are homeomorphic.
The proof is, as described above, a sequence of triples $(\sigma,i,b)$, where $\sigma$ is an isomorphism signature, $i$ is a number indicating a boundary edge, and $b$ is a bit indicating whether to fold along that boundary edge or layer along it.
This allows one to promote, as above, a polynomial-sized certificate of a property of $\TT'$ to a polynomial-sized certificate of a property of $\TT$.

As simple as the above algorithm is, it has this disadvantage, that it increases the number of tetrahedra.
If the given triangulation is assumed to be irreducible and $\partial$-irreducible, then one can do much better via crushing to a 0-efficient triangulation.
The following are essentially restatements and summaries of results from \cite{JR03}.

\begin{lemma}\label{lem:0FVTX}
  Suppose $\TT$ is a compact material triangulation that is not 0-efficient.
  Then $\TT$ admits a fundamental surface $\Sigma$ that is either an embedded $P^2$, a non-vertex-linking $S^2$, or a non-vertex-linking $D^2$.
\end{lemma}


\begin{proof}
  Suppose $\TT$ is a compact material triangulation that is not 0-efficient.
  Then $\TT$ admits a connected non-vertex-linking normal surface $\Sigma$ with $\chi(\Sigma) > 0$.
  Suppose $\Sigma = X + Y$ for nonempty normal surfaces $X,Y \neq 0$.
  Then $\chi(\Sigma) = \chi(X) + \chi(Y) > 0$.
  Hence either $\chi(X) > 0$ or $\chi(Y) > 0$.
  Thus some component $k$ of $X$ or $Y$ has $\chi(k) > 0$.
  Now, since $\Sigma$ is connected and non-vertex-linking, no component of $X$ or $Y$ is a vertex link.
  Hence in particular, $k$ is not a vertex link.
  However, since $X$ and $Y$ are nonempty, $|k \cap \TT^{(1)}| < |\Sigma \cap \TT^{(1)}|$.
  Thus by descent on $|\Sigma \cap \TT^{(1)}|$ (the \emph{weight} of $\Sigma$), $\TT$ admits a connected non-vertex-linking fundamental surface $\Sigma$ with $\chi(\Sigma) > 0$.
  The lemma follows.
  \qed
\end{proof}

The existence of the 0-efficient triangulation in the next proposition can be established by \cite[Proposition 5.9, Theorem 5.17]{JR03}; the certificate is essentially due to Schleimer (see \cite[Theorem 16.1]{Schleimer11}).

\begin{proposition}\label{prp:0Fcert}
  Suppose $\TT$ is a triangulation of a compact, orientable, $P^2$-irreducible, $\partial$-irreducible 3-manifold.

  Then there is a 0-efficient triangulation $\TT'$ of the same 3-manifold with $|\TT'| \leq |\TT|$, and also there is a certificate that $\TT \homeo \TT'$ verifiable in polynomial time.
\end{proposition}

Note carefully that we do not claim there is a polynomial-size certificate that the 0-efficient triangulation $\TT'$ is in fact 0-efficient---we only provide a certificate that $\TT \homeo \TT'$.

This is also not the most computationally efficient way to go about the certification.
Instead we have tried to give a simple certification.

\begin{proof}
  Suppose $\TT$ is a triangulation of a compact, orientable, $P^2$-irreducible, $\partial$-irreducible 3-manifold.

  Define a sequence of triangulations $\TT_0, \cdots, \TT_k$ and $k < |\TT|$ homeomorphic to $\TT$ as follows.
  Start with $i = 0$ and $\TT_i = \TT$.
  When we find $\TT_i$ is 0-efficient, we set $k = i$ and stop.
  Otherwise, by Lemma \ref{lem:0FVTX} there is a fundamental nontrivial normal surface $\Sigma_i$ with $\chi(\Sigma_i) > 0$.
  Since $\TT_i$ is $P^2$-irreducible, $\Sigma_i$ is not an embedded projective plane.
  So $\Sigma_i$ is either a sphere or a disc.
  The proof of Theorem \ref{thm:burton_crush} shows that in fact the crushing $\crush \Sigma_i$ is not just gotten from $\TT_i$ by the given collection of operations, but from $\TT_i$ surgered along $\Sigma_i$.
  So let $\TT_{i+1}'$ be $\TT_i$ surgered along $\Sigma_i$.
  Since $\TT_i$ is $P^2$-irreducible, some component $T_{i+1}$ of $\TT_{i+1}'$ is homeomorphic to $\TT_i$.
  Let $S_{i+1}$ be the other component.
  Then $S_{i+1}$ is either $S^3$ or $D^3$ according as $\Sigma_i$ is $S^2$ or $D^2$.
  Use Proposition \ref{prp:exterior} to triangulate the exterior of $\Sigma_i$, then cone off the appropriate sphere boundary component to triangulate $S_{i+1}$.
  Let $C_{i+1}$ be a certificate from Theorem \ref{thm:s3np} that $S_{i+1}$ is $S^3$ or that the double of $S_{i+1}$ is $S^3$, according as $\Sigma_i$ is $S^2$ or $D^2$.
  Next, let $R_{i+1}$ be the part of the crushing $\crush \Sigma_i$ coming from applying the operations of Theorem \ref{thm:burton_crush} to $T_{i+1}$.
  Then by the irreducibility conditions on $\TT_i$ and thus $T_{i+1}$, and by Theorem \ref{thm:burton_crush}, either $R_{i+1}$ is homeomorphic to $\TT_{i+1}$ or $R_{i+1}$ is empty.
  In the former case, let $\TT_{i+1} = R_{i+1}$.
  In the latter case, $R_{i+1}$ is either $S^3$, $D^3$, or $L(3,1)$;
  we set $\TT_{i+1}$, respectively, to the triangulation with isomorphism signature either \texttt{bkaagb}, \texttt{bGab}, or \texttt{cMcabbjak}.
  (The first two have one tetrahedron; the third has two.
  All three are 0-efficient.)
  
  By Theorem \ref{thm:s3np} and Lemma \ref{lem:hlp}, the size of each pair $(T_i,\Sigma_i,C_i)$ is bounded above by $p(|\TT_i|)$ for some polynomial $p$.
  We may assume $p$ is increasing.
  Crushing a triangulation along a nontrivial normal surface yields a triangulation with strictly fewer tetrahedra.
  Hence $p(|\TT_i|) \leq p(|\TT|)$.
  So the size of the sequence of these pairs is at most $|\TT|\cdot p(|\TT|)$, a polynomial in $|\TT|$.

  To verify that the sequence establishes a homeomorphism from $\TT$ to $\TT_k$, it suffices to verify that the sequence establishes a homeomorphism from $\TT_i$ to $\TT_{i+1}$ for each $0\leq i < k$.
  To that end, first we verify that $\Sigma_i$ is a sphere or disc in polynomial time using Proposition \ref{prp:aht};
  then use Proposition \ref{prp:exterior} as above to surger along $\Sigma_i$
  then use $C_i$ to verify in polynomial time that the appropriate component of the surgery is $S^3$ or $D^3$.
  This establishes that the other component of the surgery (not yet constructed) is homeomorphic to $\TT_i$.
  Next determine the part $R_{i+1}$ of the crushing corresponding to this other component in polynomial time.
  If $R_{i+1}$ is nonempty, verify in polynomial time that $R_{i+1}$ is isomorphic to $\TT_{i+1}$.
  Otherwise, let $K(M)$ be the predicate on 3-manifolds $M$ that ``$M$ is gotten from $\emptyset$ by finitely many operations of Theorem \ref{thm:burton_crush}.''
  By Theorem \ref{thm:burton_crush} we have proven $K(R_{i+1})$, and therefore $K(\TT_i)$.
  If $\TT_{i+1}$ has isomorphism signature \texttt{bkaagb}, which we can test in polynomial time, then we verify that $\TT_i$ is closed and a homology $S^3$.
  The only such manifold satisfying $K$ is $S^3$ itself.
  Since $\TT_{i+1}$ is homeomorphic to $S^3$, $\TT_i$ is homeomorphic to $\TT_{i+1}$.
  Similar arguments hold for \texttt{bGab} and \texttt{cMcabbjak}.
  \qed
\end{proof}

\subsection{Closed essential surface certificates}

To certify a manifold is toroidal, we need to give a certificate that a torus or Klein bottle is essential.
In an earlier version of this paper, we gave examples of such certificates.
We also conjectured the following.
\begin{theorem}\label{thm:handlebody}
  {\sc Handlebody recognition} is in \coNP{}\ among irreducible 3-manifolds, and {\sc Handlebody recognition} is in \coNP{}\ modulo GRH.
\end{theorem}
After we posted this version, Marc Lackenby kindly pointed out the following theorems from the latest revision of \cite{Lackenby16} can be used to prove Theorem \ref{thm:handlebody}:
\begin{theorem}[{\cite[Thm. 1.5]{Lackenby16}}]\label{thm:lknby_thurston}
\textsc{Thurston norm of a homology class} is in \NP{}.
\end{theorem}
\begin{theorem}[{\cite[Thm. 1.6]{Lackenby16}}]\label{thm:irred_np_lknby}
  The decision problem \textsc{irreducibility of a compact orientable 3-manifold with toroidal boundary and $b_1 > 0$} is in \NP{}.
\end{theorem}
\begin{theorem}[{\cite[Thm. 1.7]{Lackenby16}}]\label{thm:incomp_lknby}
  \textsc{Incompressible boundary} is in \NP{} $\cap$ \coNP{}.
\end{theorem}
\begin{definition}
Recall that, for any compact oriented 3-manifold $M$, the function $\chi_-: H_2(M,\partial M; \Z) \to \mathbb{Z}$ defined as
\[\chi_-([F]) = \min_{F \in [F]}\{\chi_-(F)\}\]
is called the \emph{Thurston norm},
where the minimum is taken over all oriented properly embedded surfaces $F$ in the homology class $[F]$,
where $\chi_-(S) = \max(0, -\chi(S))$ for connected surfaces $S$,
and where in general $\chi_-(F)$ is the sum of $\chi_-(S)$ over all components $S$ of $F$.

In general $\chi_-$ is only a pseudonorm, but we still call it the Thurston norm.
\end{definition}

We prove Theorem \ref{thm:handlebody} below.
We then found that it is straightforward to use it and Lackenby's new results to prove Theorems \ref{thm:ibundle} and \ref{thm:ibundle_conp}, and Corollary \ref{cor:ess_surf}. The following theorem is included both for the sake of completeness and to use a framing device to discuss these other two results. 
(We thank Nathan Dunfield for a conversation in which he brought the problem of $I$-bundle recognition to our attention.)

\begin{theorem}\label{thm:ibundle}
  For every compact surface $\Sigma$, {\sc $\Sigma\times I$ recognition} and {\sc $\Sigma\,\tilde{\times}\,I$ recognition} are in {\NP{}}.
\end{theorem}

\begin{theorem}\label{thm:ibundle_conp}
 For every compact surface $\Sigma$, {\sc $\Sigma\times I$ recognition} and {\sc $\Sigma\,\tilde{\times}\,I$ recognition} are in {\coNP{}} among orientable irreducible 3-manifolds.
\end{theorem}

\begin{corollary}\label{cor:ess_surf}
  {\sc Closed essential surface recognition} is in \NP{}\ $\cap$ \coNP{}\ for normal surfaces $\Sigma$ in connected, orientable, irreducible, $\partial$-irreducible 3-manifolds.
\end{corollary}

Note carefully that in Theorems \ref{thm:ibundle} and \ref{thm:ibundle_conp} and Corollary \ref{cor:ess_surf}, the size of the certificate returned by the algorithm is polynomial in the input triangulation and in $|\chi(\Sigma)|$, and in the case of Corollary \ref{cor:ess_surf} also polynomial in the bit-size of $\Sigma$, i.e. $\log wt(\Sigma)$.
It may be possible to make Corollary \ref{cor:ess_surf} depend polynomially only on the bit-size, but we do not pursue that here.

The careful reader will note that the \NP{}\ portion of Theorem \ref{thm:ibundle} follows almost immediately from the normal-almost-normal sweepout certificates in section 6 of Schleimer's 2001 PhD thesis \cite{SchleimerPhD}.
Moreover, it is relatively easy to show from Theorem 3 of \cite{Ivanov08} that one can extend the $\TorusxI$ recognition algorithm of \cite{Haraway14} to show that {\sc $\TorusxI$ recognition} is in \NP{}.
We regard these approaches as more viable for implementation with current software libraries like Regina (see \cite{Regina}).
However, the result also follows directly from Lackenby's work, so for ease of reference we rely on \cite{Lackenby16}.

We turn to the proofs of these results now.

\begin{proof}[Theorem \ref{thm:handlebody}]
  Let $\TT$ be a triangulation with boundary $F$ a genus $g$ surface.
  Suppose $\TT$ is not a handlebody of genus $g$.

  If $H_1(\TT,\Z) \not\cong \Z^g$, then $\TT$ is not a handlebody.
  Since homology is computable in polynomial time, $\TT$ may be verified not to be a handlebody in polynomial time, with no certificate necessary.

  Otherwise, we start out assuming $\TT$ is irreducible and will address the reducible case afterward.

  In the irreducible case, find a maximal set $\hat{D}$ of non-isotopic normal (\emph{sic}) compressing disks in $\TT$.
  By the work of Jaco and Tollefson, specifically \cite[Theorem 6.2]{JT95}, we may arrange that the disks of $\hat{D}$ are all fundamental surfaces.
  (This is because Jaco and Tollefson construct vertex-normal disks.
  Briefly, a vertex-normal disk $S$ is a normal disk admitting no sum of the form $S + \cdots + S = X + Y$ with $X,Y$ not consisting of disjoint copies of $S$.
  Vertex-normal disks are always fundamental.)
  So by Lemma \ref{lem:hlp} and Haken finiteness, the union of these disks has polynomial bit-size.
  By Proposition \ref{prp:exterior} above, construct on the exterior $\TT - \hat{D}$ a triangulation $\hat{\TT}$ of size polynomial in $|\TT|$.

  Suppose every component of $\hat{\TT}$ had $\beta_1 = 0$.
  Then each component would have a sphere boundary component.
  By irreducibility, each component would be a ball, and $\TT$ would be a handlebody, contrary to assumption.
  So some component $K$ of $\hat{\TT}$ has $\beta_1 > 0$, and is irreducible.
  If $\partial K$ were compressible, then it would admit a compressing disk.
  This would imply that $\TT$ would admit yet another normal disk disjoint from $\hat{D}$ and not normally isotopic to any element of $\hat{D}$, contrary to the maximality of $\hat{D}$.
  Thus $\partial K$ is incompressible.
  We may construct a certificate $c_{inc}$ of this fact of size polynomial in $|\TT|$ by Theorem \ref{thm:incomp_lknby}.
  Together, we contend $\hat{D}$, $\hat{\TT}$, $K$, and $c_{inc}$ constitute a certificate that $\TT$ is not a handlebody.

  Finally, if $\TT$ is reducible, then we run the standard argument of constructing a polynomially sized reduction along a reducing sphere using Proposition \ref{prp:exterior} above, coning off the sphere's sides, and applying \cite[Theorem 11.2]{Zentner16}, which depends upon GRH.

  To verify an empty certificate, just determine whether or not $\TT$ is a homology handlebody.
  If not, the certificate is verified; otherwise, the verification fails.
  To verify the certificate constructed when $\TT$ is an irreducible homology handlebody, verify that $\TT$ surgers along $\hat{D}$ into $\hat{\TT}$; then verify that $K$ is a component of $\hat{\TT}$; then verify that $c_{inc}$ is a certificate that $\partial K$ is incompressible.
  All these verifications may be performed in time polynomial in the sizes of the certificates.
  By \cite{AHT06} we may verify as well that $\hat{D}$ is a disjoint union of discs.
  Thus we may verify that $K$ is gotten from $\TT$ by surgering along discs.
  But we have also verified that $K$ has incompressible boundary.
  Furthermore, we may calculate $\beta_1(K)$; if it is 0, the verification fails.
  Otherwise, $\beta_1(K) > 0$.
  There is no handlebody with $\beta_1 > 0$ and with incompressible boundary.
  Hence $K$ is not a handlebody.
  However, surgering a handlebody along a disc always produces another handlebody (or two).
  Hence $\TT$ cannot be a handlebody.
  Thus the given certificate demonstrates $\TT$ is not a handlebody.

  Finally, to verify the reducibility certificate, verify that the given surface is a sphere, that the reduction along that sphere is as given in the certificate, and that the reduction is nontrivial using the verification of \cite[Theorem 11.2]{Zentner16}.
  \qed
\end{proof}

\begin{proof}[Theorem \ref{thm:ibundle}]
  Suppose $\Sigma$ is a compact surface.

  First, suppose $\TT$ triangulates a manifold of the form $\Sigma \times I$, with $\Sigma$ orientable.
  If $\Sigma$ has nonempty boundary, then $\TT$ is a handlebody with the same Euler characteristic as $\Sigma$.
  Theorem 3 of \cite{Ivanov08} provides a certificate $h_\TT$ that $\TT$ is a handlebody.
  To verify that $\TT$ triangulates $\Sigma\times I$, verify with $h_\TT$ that $\TT$ is a handlebody, then check, in polynomial time, that $\TT$ has the same Euler characteristic as $\Sigma$.
  On the other hand, if $\Sigma$ is closed, then $(\TT,\emptyset)$ is a product sutured manifold.
  Let $\HH$ be the handle structure dual to $\TT$, i.e.~the handle structure gotten from thickening the dual spine of $\TT$.
  One may construct this handle structure in time polynomial in $|\TT|$.
  Then $\HH$ is of uniform type, as required in Theorem 12.1 of \cite{Lackenby16}, since its 0-handles' intersection patterns with the 1- and 2-handles are subtetrahedral (in fact, tetrahedral, being the duals of the tetrahedra of the triangulation).
  So Theorem 12.1 of \cite{Lackenby16} applies, and $(\HH, \emptyset)$ admits a polynomial-time verifiable certificate $ps_\TT$ that it is a product sutured manifold.
  To verify that $\TT$ triangulates $\Sigma\times I$, verify with $ps_\TT$ that it is a product sutured manifold, then check, in polynomial time, that its two boundary components are homeomorphic to $\Sigma$.

  Second, suppose that $\TT$ triangulates a twisted $I$-bundle $\Sigma\tilde{\times}I$ with $\Sigma$ nonorientable.
  Again, if $\Sigma$ has nonempty boundary, one may verify the product as before using Ivanov's handlebody certificate.
  If instead $\Sigma$ is closed, then $\TT$ admits a double cover $\tilde{\TT}\homeo \tilde{\Sigma}\times I$ that is a product $I$-bundle over the orientation cover of $\Sigma$.
  One may certify that $\tilde{\TT}$ is $\tilde{\Sigma}\times I$, and this constitutes a proof that $\TT$ is $\Sigma\tilde{\times}I$.
  Thus, {\sc $\Sigma\times I$ recognition} and {\sc $\Sigma\tilde{\times} I$ recognition} are in \NP{}\ among orientable 3-manifolds.
\end{proof}

Before getting to the proof of Theorem \ref{thm:ibundle_conp} we give the following characterizing lemma. 
  
\begin{lemma}\label{lem:bundle_norms}
  Let $\Sigma$ be a closed orientable surface.
  Let $M$ be an irreducible orientable 3-manifold with integer homology equivalent to $\Sigma \times I$.
  
  $M$ is homeomorphic to $\Sigma \times I$ if and only if $\chi_- = 0$ on all $H_2(M,\partial M;\Z)$.
\end{lemma}

\begin{proof}
  Suppose an oriented irreducible 3-manifold $M$ is an integral homology $\Sigma \times I$.
  
  For the only-if-direction, suppose $M$ is homeomorphic to $\Sigma \times I$.
  Then every element of $H_2(M,\partial M; \Z)$ is representable with a disjoint union of annuli, which has Thurston norm 0.
  So $\chi_- = 0$ on $H_2(M,\partial M; \Z)$.

  For the if-direction, suppose $\chi_- = 0$ on all $H_2(M, \partial M; \Z)$.
  First recall the long exact sequence in homology for the pair $(M, \partial M)$.
  We have that $j_\ast:H_2(M;\Z) \to H_2(M,\partial M; \Z)$ is identically 0, since $M$ is an integer homology $\Sigma \times I$.
  So closed surfaces are nullhomologous in $H_2(M,\partial M; \Z)$.
  Moreover, $\partial_2: H_2(M, \partial M; \Z) \to H_1(\partial M; \Z)$ is injective. 
  Likewise, $j_\ast: H_1(M;\Z) \to H_1(M,\partial M; \Z)$ is also identically 0.
  Thus $H_2(M,\partial M; \Z) \simeq \image \partial_2 = \ker i_{\partial \ast}$,
  where $i_{\partial \ast}: H_1(\partial M; \Z) \to H_1(M; \Z)$.
  Also, curves in $\partial M$ nullhomologous in $M$ are nullhomologous in $\partial M$.
    
  Now, let $B$ be a basis for $H_2(M,\partial M; \Z)$.
  Suppose $b \in B$.
  Then $\chi_-(b) = 0$.
  Let $A_b$ be an oriented properly embedded surface with
  \begin{equation}\label{eq:ab}
    [A_b] = b, \qquad \chi_-(A_b) = 0.
  \end{equation}

  All the components of $A_b$ are oriented and hence orientable.
  Any closed components of $A_b$ are nullhomologous in $H_2(M, \partial M; \Z)$, and are spheres or tori.
  Omitting them preserves (\ref{eq:ab}).
  The boundary of any disc in $A_b$ is nullhomologous in $M$, and therefore nullhomologous in $\partial M$.
  Thus any disc in $A_b$ is nullhomologous in $H_2(M,\partial M; \Z)$.
  Thus omitting discs also preserves (\ref{eq:ab}).
  So we may assume $A_b$ is a disjoint union of annuli.

  Suppose $A$ is an oriented recurrent annulus properly embedded in $M$, i.e.~an annulus with $\partial A$ lying all on one boundary component.
  Then $i_\ast(\partial_2([A])) = 0$.
  Let $\Sigma_\pm$ be the boundary components of $M$.
  On each of the subspaces $H_1(\Sigma_\pm; \Z)$, the map $i_\ast$ restricts to an isomorphism.
  Therefore $\partial_2([A]) = 0$.
  Thus $\partial A$ bounds an oriented surface in $\partial M$.
  Attaching that surface to $A$ yields a closed surface homologous to $A$ in $H_2(M,\partial M;\Z)$.
  So $[A] = 0$ in $H_2(M,\partial M;\Z)$.
  Thus we may omit recurrent annuli from $A_b$ preserving (\ref{eq:ab}).

  Finally, suppose $\gamma$ is a boundary component of a non-recurrent annulus $A$ properly embedded in $M$ such that $\gamma$ is nullhomologous in $M$.
  Then $\gamma$ is nullhomologous in $\partial M$.
  Therefore $\gamma$ bounds an oriented surface $Y$ in $\partial M$ such that $Y \cup A$ is an oriented surface homologous to $A$ in $H_2(M, \partial M; \Z)$.
  Pushing this surface off from that boundary component of $M$ yields a surface bounded by the other boundary curve of $A$.
  The same argument now applies to that boundary curve, so that $A$ is homologous to a closed surface in $H_2(M,\partial M; \Z)$.
  Thus $A$ is nullhomologous.
  Therefore we may omit any annulus with a nullhomologous boundary curve from $A_b$ preserving (\ref{eq:ab}).

  Let $\mathbf{A}$ be the set of all the oriented annulus components of all the $A_b$.
  Let $[\mathbf{A}]$ be the set of all their homology classes.
  Then the span of $[\mathbf{A}]$ is $H_2(M,\partial M; \Z)$ since $B$ is a basis.
  Select a new basis from $[\mathbf{A}]$, and let $\mathbf{A}'$ be the set of the corresponding annuli.
  We may isotope these annuli so that they are all transverse.
  By a standard argument, after isotoping the annuli to lie in minimal position, they only intersect in core curves or essential arcs.

  Suppose two distinct annuli $A,A'$ intersect in a core curve.
  One resolution of this intersection curve yields two recurrent annuli.
  That means $[\partial A] = \pm[\partial A']$.
  That is, $\partial_2([A]) = \pm \partial_2([A'])$.
  Since $\partial_2$ is injective, that means $[A] = \pm [A']$.
  This contradicts $[\mathbf{A}]$ being a basis.
  So no two distinct annuli intersect in core curves.
  Hence, after isotoping the annuli to lie in minimal position, they intersect only in essential arcs.
  
  Now, the union of their boundary curves on $\partial M$ is filling.
  That is, its complement consists of discs.
  Let $X$ be a regular neighborhood of the union of $\partial M$ with these annuli.
  The boundary of $X$ consists of $\partial M$ and spheres.
  By irreducibility the spheres must bound balls, and they must bound them away from $\partial M$.
  Hence $M$ is $\Sigma \times I$.

  \qed
\end{proof}

\begin{proof}[Theorem \ref{thm:ibundle_conp}]
  Suppose first that $\Sigma$ is orientable.
  
  Suppose $\TT$ is a triangulation of a compact orientable irreducible 3-manifold that is not $\Sigma\times I$.
  
  If $\TT$ is not an integer homology $\Sigma\times I$, we can use the empty certificate.

  Otherwise $\TT$ is an integer homology $\Sigma\times I$.
  Let $\Sigma_+$ be one of the boundary components of $M$.
  Select a basis for $H_1(\Sigma_+;\Z/2\Z)$ every element of which is represented by a simple closed normal curve intersecting every edge of $\Sigma_+$ at most once.
  Let $\Gamma$ be a set consisting of one such curve from each homology class in the basis.
  Picking some orientation of these curves, we get a basis $B_+$ for $H_1(\Sigma_+; \Z)$ of size polynomial in $\TT$ and $|\chi(\Sigma)|$.

  Let $i_\pm: \Sigma_\pm \to M$ be the inclusions.
  By assumption, these induce isomorphisms on homology.
  Define the map $k: H_1(\Sigma_+; \Z) \to H_1(\partial M; \Z)$ by
  $k(c) = c - i_{-\ast}^{-1}(i_{+\ast}(c))$.
  Then $B' = \{k(b_+)\ |\ b_+ \in B_+\}$ is a basis for $\ker i_\ast$, where $i_\ast: H_1(\partial M;\Z) \to H_1(M;\Z)$.
  However, by the long exact sequence of the pair $(M,\partial M)$, since $M$ is a homology $\Sigma \times I$, $\ker i_\ast = \image \partial_2$.
  So $B'' = \partial_2^{-1}(B')$ is a basis for $H_2(M,\partial M; \Z)$.
  Since $B_+$ had size polynomial in $\TT$ and $|\chi(\Sigma)|$, and since the above maps on homology and their inverses are representable with matrices with polynomially many coefficients of polynomial size, likewise $B''$ has size polynomial in $\TT$ and $|\chi(\Sigma)|$.

  Since $\TT$ is an orientable irreducible integer homology $\Sigma \times I$ but isn't $\Sigma \times I$, by the if-direction of Lemma \ref{lem:bundle_norms}, $\chi_-$ is not identically 0 on $H_2(M,\partial M; \Z)$.
  Therefore, $\chi_-$ is not identically 0 on $B''$.
  So we may pick an element $b \in B''$ with $\chi_-(b) > 0$.
  A natural isomorphism $L: H_2(M,\partial M; \Z) \to H^1(M;\Z)$ representable by a matrix of polynomially many coefficients of polynomial size is given by Lefschetz duality.
  (Theorem \ref{thm:lknby_thurston} is phrased in terms of $H^1(M;\Z)$, not $H_2(M,\partial M;\Z)$.)
  So $L(b)$ is of size polynomial in $\TT$ and $b$, i.e.~in $\TT$ and $|\chi(\Sigma)|$.
  By Theorem \ref{thm:lknby_thurston}, there is a certificate $C$ verifiable in time polynomial in $\TT$ and $L(b)$, i.e.~in $\TT$ and $\chi(\Sigma)$, that $\chi_-(b) > 0$.

  To verify, suppose one is given $\TT$, and either an empty certificate, or a simplicial 1-cocycle $\eta$ and a certificate $C$ that $\chi_-(L(\eta)) > 0$.
  If the certificate is empty, check in time polynomial in $\TT$ and $|\chi(\Sigma)|$ that the integral homology of $\TT$ differs from $\Sigma \times I$.
  If that is verified, then $\TT$ is not $\Sigma \times I$.
  Otherwise, use $C$ to check in time polynomial in $\TT$ and $|\chi(\Sigma)|$ that $\chi_-(L(\eta)) > 0$.
  If that is verified, then by the only-if direction of Lemma \ref{lem:bundle_norms}, $\TT$ is not $\Sigma \times I$.

  Finally, if $\Sigma$ is nonorientable, we may verify that $\TT$ is not $\Sigma\tilde{\times} I$ by showing either that $\TT$ does not admit a double cover, or by showing that its unique connected double cover is not $\Sigma\times I$.
\qed
\end{proof}

We remind the reader that the verification in the above proof does not assume $\TT$ is irreducible.
Only the construction of the certificate makes that assumption.

\begin{proof}[Corollary \ref{cor:ess_surf}]
  Suppose $\Sigma$ is a closed essential normal surface in $\TT$ with $\TT$ orientable, irreducible, and $\partial$-irreducible.
  That is, suppose $\Sigma$ is incompressible and not boundary parallel.
  ($\Sigma$ is automatically $\partial$-incompressible, being closed, and this is why we restrict to the closed case.)
  If $\Sigma$ is nonorientable, its normal double is its orientation cover, and is essential if and only if $\Sigma$ is.
  Thus we may assume $\Sigma$ is orientable.
  To certify $\Sigma$ is incompressible, take a polynomial-time verifiable triangulation $\TT'$ of its exterior with Proposition \ref{prp:exterior}, and use Theorem \ref{thm:incomp_lknby} to verify that $\TT'$ has incompressible boundary.
  To show $\Sigma$ is not boundary-parallel, either show that $\TT'$ is connected, or for both components of $\TT'$ use Theorem \ref{thm:ibundle_conp} to verify that they are not $\Sigma \times I$.
  The latter certificate will be verifiable in time polynomial in $\TT'$ and $\Sigma$, which is polynomial in $\TT$ and $\Sigma$.
  \qed
\end{proof}

\subsection{Toroidality}
We now consider the other decision problems discussed in the introduction.
We begin with Theorem \ref{main_thm_toroidal}, on toroidality.
First we require the following proposition.

\begin{proposition}\label{prp:tornp}
  {\sc Toroidal recognition} lies in \NP{}\ among compact, connected, irreducible, $\partial$-irreducible, orientable 3-manifolds admitting no essential Klein bottles.
\end{proposition}

\begin{proof}
  Suppose $\TT$ is a triangulation of a compact, connected, irreducible, $\partial$-irreducible, orientable 3-manifold that has no Klein bottle and is toroidal.
  By Proposition \ref{prp:fundFaults}, there is a fundamental essential torus $F$ in $\TT$.
  By Lemma \ref{lem:hlp}, $F$ is representable in size polynomial in $\TT$.
  By Proposition \ref{prp:exterior}, we may construct a triangulation $\TT'$ of the exterior of $F$ in time polynomial in the sizes of $\TT$ and $F$, hence in time polynomial in $\TT$.

  If $\TT'$ is connected, then $\TT'$ is $\partial$-irreducible, since $\TT$ is $\partial$-irreducible and $F$ is essential.
  Lackenby's Theorem \ref{thm:incomp_lknby} provides a certificate that $\TT'$ is $\partial$-irreducible.
  Given such a certificate, we can verify that $F$ is incompressible in time polynomial in $\TT$, and hence that $\TT$ is toroidal.
  
  Otherwise, since $F$ is essential and $\TT$ is connected, $\TT'$ has two components, $L$ and $R$, neither of which is homeomorphic to $\TorusxI$, neither of which is homeomorphic to $\SxD$, and both of which are irreducible (and hence $\partial$-irreducible, not being $\SxD$).
  The certificate we return is the quintuple $(F, L_T, R_T, L_D, R_D)$, where $L_T,R_T$ are the certificates guaranteed by Theorem \ref{thm:ibundle_conp} that $L$ and $R$ are not $\TorusxI$, and where $L_D,R_D$ are the certificates guaranteed by Theorems \ref{thm:irred_np_lknby} and \ref{thm:incomp_lknby} that $L,R$ are irreducible and are not solid tori.
  By Theorem \ref{thm:ibundle_conp} and Lackenby's Theorems \ref{thm:irred_np_lknby} and \ref{thm:incomp_lknby},
  verifying these certificates takes time polynomial in $\TT'$, and hence polynomial in $\TT$.
  This verifies that $\TT$ is toroidal.
  \qed
\end{proof}

We may now prove a theorem asserted in the introduction.


\begin{proof}[Theorem \ref{main_thm_toroidal}]
  Suppose $\TT$ is a triangulation of a connected, orientable, irreducible, $\partial$-irreducible, toroidal 3-manifold.
  Either $\TT$ contains an essential Klein bottle or not.
  If it does, then by Proposition \ref{prp:fundFaults} there is a fundamental one.
  Corollary \ref{cor:ess_surf} produces a certificate verifiable in time polynomial in $\TT$ that this Klein bottle is essential.
  Otherwise, since $\TT$ is toroidal, by Proposition \ref{prp:tornp} there is a certificate verifiable in time polynomial in $\TT$ that $\TT$ contains an essential torus.
  In the latter case we already have a proof of toroidality.
  In the former case, the double of the Klein bottle is a torus.
  Since the Klein bottle is essential, the torus is necessarily essential.
  Since the Klein bottle is fundamental, the torus has size polynomial in $\TT$.
  Either option is a certificate of toroidality.

  Since $\TT$ admits an essential torus, $b_1(\TT) > 0$.
  So a certificate of irreducibility follows immediately from Lackenby's work, stated here as Theorem \ref{thm:irred_np_lknby}.
  \qed
\end{proof}

Corollary \ref{thm:satellite} follows immediately; it only requires going from a knot diagram to a triangulation, which is a standard argument of Hass, Lagarias, and Pippenger \cite[Lemma 7.1]{HLP99}.
This result was achieved by Baldwin and Sivek in \cite{BaldwinSivek17} under the additional assumption of the Generalized Riemann Hypothesis.

\subsection{Essential \Mob s and annuli}\label{sec:annuli}

We now discuss the final situations: when the manifold contains an essential \Mob\ or annulus.


\begin{proposition}\label{prp:m2np}
  Suppose a compact, orientable, $P^2$-irreducible, $\partial$-irreducible triangulation $\TT$ has no essential tori or Klein bottles.
  If $\TT$ admits a properly embedded \Mob, then there is a certificate of this fact verifiable in time polynomial in $\TT$.
\end{proposition}

\begin{proof}
  Suppose $\TT$ admits an embedded \Mob\ $M$.

  If $M$ is compressible, then $\TT$ admits an embedded $P^2$, contrary to $P^2$-irreducibility.
  If $M$ is incompressible but is $\partial$-compressible, let $D$ be a $\partial$-compressing disk for $M$.
  Then surgering $M$ along $D$ yields another disk $Z$.
  Since $M$ is incompressible, $Z$ is a compressing disk for $\partial \TT$, contrary to $\partial$-irreducibility.
  Thus $M$ must be incompressible and $\partial$-irreducible.
  Therefore $M$ is not $\partial$-parallel, and hence $M$ is essential.
  By Proposition \ref{prp:fundFaults} there is a fundamental normal embedded \Mob\ $\mu$ in $\TT$.
  By Lemma \ref{lem:hlp} and Proposition \ref{prp:aht}, $\mu$ constitutes a certificate verifiable in time polynomial in $\TT$ that $\TT$ admits an embedded \Mob.
  \qed
\end{proof}

The following is the main result of this subsection.

\begin{proposition}\label{prp:annulusxNP}
  Suppose $\TT$ is compact and orientable.

  Suppose further that $\TT$ has no essential surfaces of nonnegative Euler characteristic apart from annuli.

 If $\TT$ admits an essential embedded annulus, there is a certificate verifiable in time polynomial in $\TT$ of this fact. 
\end{proposition}


We prove this at the conclusion of this subsection, using the following three results.


\begin{lemma}\label{lem:nonsepAs}
  Suppose a 0-efficient triangulation $\TT$ admits a non-separating annulus $A$.
  Then $A$ is essential.
\end{lemma}


  This holds simply because inessential annuli are separating in irreducible, $\partial$-irreducible 3-manifolds.

We now show that if a non-separating annulus exists in the manifold, there is one of manageable size for our purposes.


\begin{lemma}\label{lem:nonsepAx}
  Suppose $\TT$ is compact, orientable, and 0-efficient.
  Suppose further that $\TT$ has no essential surfaces of nonnegative Euler characteristic apart from annuli.
  Finally, suppose $\TT$ admits a non-separating annulus.

  Then $\TT$ admits a fundamental non-separating annulus.
\end{lemma}

\begin{proof}
  Suppose $A$ is a non-separating annulus in $\TT$.
  Being essential by Lemma \ref{lem:nonsepAs}, $A$ must admit an isotopic normal representative, the isotopy coming from the normalization procedure of \cite{JR03} and the fact that $\TT$ is irreducible and $\partial$-irreducible by assumption. 
  Suppose that $A = X + Y$ for nontrivial normal surfaces $X,Y$, such that $X \cap Y$ has a minimal number of components.
  By a standard trick (see for example \cite[Lem.~3.3.30]{Matveev07}), we may assume $X$ and $Y$ are both connected.
  Let $[\cdot]$ denote the homology class in $H_2(\TT, \partial \TT; \mathbb{Z}/2\mathbb{Z})$.
  Then $[A] = [X+Y] = [X] + [Y]$, since Haken sum respects $\mathbb{Z}/2\mathbb{Z}$ homology.
  Since $A$ is nonseparating, $[A] \neq 0$.
  Therefore, without loss of generality, $[X] \neq 0$.
  Since $[X] \neq 0$ and $\TT$ is orientable, $X$ is nonseparating.
  Since $\TT$ is 0-efficient, $\chi(X) = \chi(Y) = 0$.
  Suppose $X$ were not an essential annulus.
  Then a minimal genus representative of $[X]$ would be a nonseparating surface of nonnegative Euler characteristic different from an annulus, contrary to assumption.
  So $X$ is an essential annulus.
  Finally, since $wt(A) = wt(X) + wt(Y)$ and $wt(Y) > 0$ by nontriviality, $wt(X) < wt(A)$.
  The result follows by descent on weight.
  \qed
\end{proof}

The following proposition handles one of the final cases we need to deal with, an essential separating annulus.

\begin{proposition}\label{prp:SxDx2}
  Suppose an orientable, irreducible, $\partial$-irreducible 3-manifold $\TT$ is obtained by identifying two solid tori along two disjoint annuli in their boundaries, and has one boundary component, a torus.
  Let $A$ be the annulus in $\TT$ gotten by identifying these annuli.

  Then $\TT$ is Seifert fibered over the disk with two exceptional fibers, and $A$ is an essential annulus fibered by regular fibers.
\end{proposition}

\begin{proof}
  Call the annuli of the solid tori $\alpha$ and $\beta$.

  Suppose first that $\alpha$ and $\beta$ are $\pi_1$-injective in their respective components.
  Then the $\DxS$ components admit Seifert fiberings such that $\alpha$ and $\beta$ are fibered by regular fibers.
  If either $\alpha$ or $\beta$ is longitudinal in its component, then $\TT$ would be $\DxS$, and hence be $\partial$-reducible, contrary to assumption.
  Therefore, neither $\alpha$ nor $\beta$ is longitudinal.
  Thus, both the Seifert fiberings have exceptional fibers.
  Therefore, $\TT$ is Seifert fibered over a disk with two exceptional fibers, and $A$ is an essential annulus fibered by regular fibers.

  Thus to conclude the proof it will suffice to assume not both $\alpha$ and $\beta$ are $\pi_1$-injective in their components, and derive a contradiction.
  As before, if but one of them were $\pi_1$-injective, then $\TT$ would have a compressible boundary component, contrary to $\partial$-irreducibility.
  Thus, suppose both $\alpha$ and $\beta$ are not $\pi_1$-injective in their respective components.
  If either were meridional, then $\TT$ would admit a non-separating sphere, contrary to irreducibility.
  Thus we may assume both $\alpha$ and $\beta$ are trivial.
  But in this case, $\TT$ would have two boundary components.
  \qed
\end{proof}

We may now conclude this subsection with a proof of its main result.
\begin{proof}[Prop.\,\ref{prp:annulusxNP}.]
  Suppose $\TT$ has no essential surfaces of nonnegative Euler characteristic, except for essential annuli.
  By Theorems \ref{thm:irred_np_lknby} and \ref{thm:incomp_lknby} there is a certificate $c_I$ that $\TT$ is a compact irreducible $\partial$-irreducible 3-manifold verifiable in time polynomial in $\TT$.

  If $\TT$ admits a nonseparating annulus, then by Lemma \ref{lem:nonsepAx}, $\TT$ admits a fundamental such annulus $A$.
  By Proposition \ref{prp:exterior}, one may construct a triangulation $\TT'$ of $\TT - A$ from $A$ and $\TT$ in time polynomial in $\TT$ and $\log |A\cap \TT^{(1)}|$.
  By Proposition \ref{prp:exterior} and Lemma \ref{lem:nonsepAs}, since $\TT'$ is connected, $c_I$ constitutes a certificate that $A$ is essential.
  Thus, $(A, c_I)$ together constitute a certificate that $\TT$ admits an essential annulus.
  Since $A$ and $c_I$ have size polynomial in $\TT$, this certificate has size polynomial in $\TT$.

  Otherwise, $\TT$ admits an essential separating annulus $A$.
  Since $A$ is separating, $\partial A \subset T$ for some boundary torus $T$ of $\TT$.
  A regular neighborhood of $A \cup T$ in $\TT$ has two torus boundary components $T', T''$.
  Since $\TT$ has no essential tori, $T',T''$ are compressible.
  They must compress away from $T$ to spheres $S,S'$.
  Since $\TT$ is irreducible, these spheres must bound balls away from $T$.
  Thus in fact $T',T''$ bound solid tori away from $T$.
  Hence $A$ separates $\TT$ into two solid tori.
  Since $A$ is essential in $\TT$, the core of $A$ is incompressible in each solid torus.
  Hence this core induces Seifert fiberings in these solid tori.
  Therefore $\TT$ is Seifert-fibered over a disk with two exceptional fibers, with $A$ being a vertical annulus.
  If $\TT$ admits a horizontal essential annulus, then $\TT$ is either $\TorusxI$ or $\KxxI$; since $\TT$ has but one boundary component, $\TT$ must be $\KxxI$.
  But we assumed $\TT$ admitted no Klein bottle, a contradiction.
  Thus essential annuli in $\TT$ are vertical.
  By Proposition \ref{prp:fundFaults} $\TT$ admits a fundamental essential annulus $A$, and by the above argument, $A$ is vertical.
  Thus $\TT - A$ is two solid tori, and $\TT$ fits the conditions of Proposition \ref{prp:SxDx2}.
  Now, by Proposition \ref{prp:exterior} we may triangulate $\TT - A$ by $\TT'$ given $A$ and $\TT$ in time polynomial in $\TT$ and $\log |A\cap \TT^{(1)}|$, and by \cite[Theorem 3]{Ivanov08}, we may construct a certificate $c_t$ that $\TT'$ is two solid tori, verifiable in time polynomial in $\TT'$, and hence polynomial in $\TT$.
  Thus, $c_I$ and $c_t$ together constitute a proof that $A$ is essential, by Proposition \ref{prp:SxDx2}.
  So $(A, c_I, c_t)$ constitutes a certificate that $\TT$ admits an essential annulus.
  Since $A$ and $c_I$ have size polynomial in $\TT$, and since $c_t$ has size polynomial in $\TT'$, and hence polynomial in $\TT$, this certificate has size polynomial in $\TT$.
  \qed
\end{proof}

\subsection{Non-hyperbolicity}
We end this section on certificates with the main theorem stated in the introduction, following the outline of Proposition \ref{prp:fundFaults}.

\begin{proof}[Thm.\,\ref{thm:main}.]
  Suppose $\TT$ is a compact orientable triangulation and $\partial \TT$ is a nonempty union of tori, and suppose $\TT$ is not hyperbolic.
  Then by Theorem \ref{thm:thurston}, $\TT$ admits an essential surface of nonnegative Euler characteristic.

  If $\TT$ admits an embedded projective plane, then by Proposition \ref{prp:fundFaults} it admits a fundamental normal such surface $\Sigma$.
  $\Sigma$ itself constitutes a proof against hyperbolicity, verifiable as an embedded projective plane in time polynomial in $\TT$ by the bounds of Lemma \ref{lem:hlp} and the algorithm of Proposition \ref{prp:aht}.

  Otherwise, if $\TT$ admits an essential sphere, then by Proposition \ref{prp:fundFaults}, since it has no embedded projective plane, it admits a fundamental normal essential sphere $\Sigma$.
  If $\Sigma$ is non-separating, then $\Sigma$ on its own constitutes a proof against hyperbolicity, verifiable in time polynomial in $\TT$ by the bounds of Lemma \ref{lem:hlp} the algorithm of Proposition \ref{prp:aht}, and the exterior algorithm of Proposition \ref{prp:exterior}, since one may determine whether or not a triangulation (in this case, the exterior of $\Sigma$) is connected in time polynomial in the triangulation.
  Otherwise, the exterior $\TT' = \TT - \Sigma$ has two components $L'$ and $R'$; cap them off with balls to get the connect-summands $L$ and $R$ of $\TT$.
  By assumption {\sc $S^3$ recognition} is in \coNP{}, so there are certificates $c_L$ and $c_R$ verifiable in time polynomial in $L$ and $R$, and hence in $\TT$, that $L$ and $R$ are not $S^3$.
  Thus $(\Sigma, c_L, c_R)$ constitutes a certificate verifiable in time polynomial in $\TT$ that $\TT$ is reducible, and hence is not hyperbolic.

  Otherwise, if $\TT$ admits a compressing disk, then by irreducibility, $\TT$ must be a solid torus.
  There is a certificate $c_d$ that $\TT$ is a solid torus verifiable in time polynomial in $\TT$ by Theorem 3 of \cite{Ivanov08}.

  Otherwise, if $\TT$ admits an essential Klein bottle or torus, then by Corollary \ref{cor:ess_surf} there is a certificate $c_{ess}$ of this fact verifiable in time polynomial in $\TT$.

  Otherwise, if $\TT$ admits an essential \Mob, then by Proposition \ref{prp:m2np}, there is a certificate $c_m$ of this fact verifiable in time polynomial in $\TT$.

  Otherwise, and finally, if $\TT$ admits an essential annulus, then by Proposition \ref{prp:annulusxNP}, either there is a certificate $c_a$ of this fact verifiable in time polynomial in $\TT$, or there is a certificate $c_a$ verifiable in time polynomial in $\TT$ that $\TT$ splits along an annulus into two solid tori.

  In all of the above cases, the given certificates show that $\TT$ is not hyperbolic.

  If $\TT$ has none of these surfaces, then since $\TT$ is Haken, having nonempty torus boundary, by a well-known result of Thurston 
  (Theorem \ref{thm:thurston})
  $\TT$ is hyperbolic contrary to assumption.
  \qed
\end{proof}\

\section{Discussion}

One initial remark is that in practice, one would want to use vertex-normal surfaces instead of fundamental normal surfaces.
Reproving the above results for vertex-normal surfaces would be a useful next task to do.

In the course of showing that {\sc $S^3$ recognition} lies in \coNP{}~modulo the Generalized Riemann Hypothesis \cite[Theorem 11.2]{Zentner16}, a good bit of the proof is allocated to constructing irreducible $\SLTC$ representations of toroidal manifold groups.
Theorem \ref{main_thm_toroidal} removes the need for splicing, and hence reduces one to the case of closed, irreducible, atoroidal 3-manifolds---that is, by Perelman's resolution to Thurston's Geometrization Conjecture, it reduces one to geometric 3-manifolds.
We can further assume that the manifold is \emph{small}, i.e.~it does not contain an incompressible surface of positive genus.
Otherwise, we can use an incompressible surface which is given by a fundamental normal surface, together with Theorem \ref{thm:incomp_lknby} as a certificate of incompressibility, that the manifold is not $S^3$.
Thus, if the following conjecture is true, then {\sc $S^3$ recognition} is in fact in \coNP{}, and the results of this paper are unconditionally true.

\begin{conjecture}\label{conj:S3inCoNP}
  If $\TT$ triangulates a small geometric integral homology sphere, and $\TT$ is not $S^3$, then there is a proof that $\TT$ is not $S^3$ of length polynomial in $\TT$.
  That is, among small geometric integral homology spheres, {\sc $S^3$ recognition} is \coNP{}.
\end{conjecture}

As mentioned above, Zentner's work specifically \cite[Theorem 11.2]{Zentner16} combines with Kuperberg's work to give a proof of this conjecture, assuming the Generalized Riemann Hypothesis.
To remove this assumption, it seems a promising line of inquiry to approach this first by considering small Seifert spaces, for the hyperbolic integral homology spheres will likely prove much more difficult to verify. 

We will conclude with a discussion of the affirmative problem.
As mentioned in the introduction, a 3-manifold $M$ is said to be hyperbolic if $M\cong \Hth/\Gamma$, such that $vol(\Hth/\Gamma) < \infty$ and $\Gamma$ is a subgroup of the isometry group acting properly discontinuously.
Now, an ideal triangulation $\TT$ of $M$ admits a \emph{strict angle structure} if it determines a map from $\RR^{3n} \to \RR^n$ that corresponds each pair of opposite edges in a tetrahedron are assigned a positive number, the dihedral angle, such that all (six) dihedral angles in a tetrahedron sum to $2\pi$ and all all dihedral angles in an edge class of $\TT$ sum to $2\pi$.
Thus, an angle structure is a solution to a system of linear equations that lies in the positive cone of the solution space.
Work of Casson and Rivin \cite{rivin1994euclidean}  (see also \cite{futer2011angled}) shows that the existence of such a structure rules out each of the obstructions to non-hyperbolicity, so that a strict angle structure on $\TT$ implies the hyperbolicity of $M$.
We point out that if a triangulation supports an angle structure, an argument analogous to Lemma \ref{lem:hlp} shows there is a fundamental solution to the angle structure equations that can be verified in polynomial time.
The relevant question is then ``Which triangulations support angle structures?''
To determine how far an input triangulation $\TT$ is from a triangulation that supports an angle structure, we ask the following question:

\begin{question}
  Is there a polynomial $P$ such that for any ideal triangulation $\TT$ of a cusped hyperbolic 3-manifold, there is an ideal triangulation $\TT'$ supporting a strict angle structure that is connected to $\TT$ by at most $P(|\TT|)$ 2-3 and 3-2 moves?
\end{question}

An affirmative answer to that question together with an affirmative answer to Conjecture \ref{conj:S3inCoNP} would provide an affirmative answer to the following:

\begin{question}
  Is the {\sc Hyperbolicity problem} for cusped manifolds in \NP{} $\cap$ \coNP{}?
\end{question}

\section*{Appendix}

\begin{appendix}
  \begin{proof}[Prop. \ref{prp:exterior}.]
    Let $\BB$ be the parallelity bundle for the handle structure on $M' = cl(\TT \setminus S)$ inherited from $\TT$ as defined above.
    Let $X$ be $cl(M' - \BB)$.
    Let $V$ be $\BB \cap X$.
    Then $V$ is the common frontier of $\BB$ and $X$ in $M'$.
    Now, $X$ is the union of the non-parallelity 0-handles of the (inverted) handle structure.
    The non-parallelity 0-handles have two isomorphism types: a truncated octahedron, and a hexagonal prism.
    $V$ is the union of the vertical 4-sided 2-cells of these non-parallelity 0-handles (i.e.~those 4-sided 2-cells not lying on $S$).
    There are at least $4\cdot t$ and at most $6\cdot t$ such 2-cells in $X$, the maximum occurring when $S$ has a normal quad in every tetrahedron.
    Therefore, $V$ can be triangulated in time linear in $t$.
    Moreover $X$ can be triangulated in time linear in $t$ agreeing with the triangulation on $V$.

    It remains to triangulate $\BB$ in time linear in $t$ and $|\chi(S)|$ in a way agreeing with $V$.
    Note that $V$ has a number of components linear in $t$.
    Since $V$ is the frontier of $\BB$ in $M'$, $\BB$ has at most as many components as $V$.
    Suppose $B$ is any component of $\BB$.
    The algorithm of Theorem \ref{thm:lknby_complement} yields for $B$ the given data from items (i)--(iii).
    The data from (iii) indicates how to identify the frontier of $B$ in $M'$ with a component of $V$.
    So it will suffice to show $B$ can be triangulated in polynomial time.

    The data from (i) and (ii) identify the homeomorphism type of $B$ in polynomial time.
    Now, $B$ is an $I$-bundle over a base surface $\Sigma_B$.
    The data from (i) and (ii) identifies the homeomorphism type of $\Sigma_B$ since $\TT$ is orientable.
    Now $\Sigma_B$ is doubly-covered by a unique subsurface $Q$ of $S$, such that $B$ is the mapping cylinder of the double-cover.
    Then $\chi(Q) = 2\cdot \chi(\Sigma_B)$.
    Let $\Theta$ be formed from $Q$ by filling in all boundary components bounding discs in $S$ disjoint from $Q$.
    We have $\chi(\Theta) \geq \chi(S)$ and $\chi(\Theta) = \chi(Q) + d$ where $d$ is the number of discs filled in.
    The number of boundary components of $Q$ is at most twice $|\pi_0(V)|$, which is linear in $t$.
    So $d$ has an upper bound linear in $t$.
    Thus $\chi(Q) = \chi(\Theta) - d \geq \chi(S) - d$ has a lower bound linear in $t$ and $\chi(S)$.
    Hence $\chi(\Sigma_B)$ has such a lower bound.
    Consequently, $|\chi(\Sigma_B)|$ has an upper bound linear in $t$ and $|\chi(S)|$.
    Since $\Sigma_B$ is connected, this means $\Sigma_B$ can be triangulated with a number of triangles linear in $t$ and $|\chi(S)|$.
    Therefore $B$ can be triangulated with a number of tetrahedra linear in $t$ and $|\chi(S)|$, a polynomial as desired.
    \qed
  \end{proof}
\end{appendix}

\bibliographystyle{plain}

\end{document}